\documentclass[12pt]{amsart}

\usepackage{amscd,amssymb}
\usepackage{amsthm,amsmath,amssymb}
\usepackage[matrix,arrow]{xy}

\newtheorem{theorem}[equation]{Theorem}

\newtheorem{lemma}[equation]{Lemma}
\newtheorem{corollary}[equation]{Corollary}

\theoremstyle{definition}

\theoremstyle{remark}
\newtheorem{remark}[equation]{Remark}

\author{Ivan Cheltsov}

\title{On nodal quintic fourfold}

\address{\begin{tabbing}
\hspace*{25 em}\=\kill
Steklov Institute of Mathematics \>School of Mathematics\\
8 Gubkin street, Moscow 117966   \>The University of Edinburgh\\
Russia                           \>Kings Buildings,  Mayfield Road\\
                                 \> Edinburgh EH9 3JZ, UK\\
\texttt{cheltsov@yahoo.com}      \>\\
                                 \>\texttt{I.Cheltsov@ed.ac.uk}
\end{tabbing}}

\thanks{The author is very grateful to I.\,Aliev, A.\,Corti, M.\,Gri\-nenko,
V.\,Is\-kov\-skikh, M.\,Mella, J.\,Park, Yu.\,Pro\-kho\-rov,
A.\,Pukhlikov and V.\,Sho\-ku\-rov for very useful and fruitful
conversations.}

\begin{document}

\begin{abstract}
We use the Shokurov connectedness principle and the Corti
inequality to prove the birational superrigidity of a nodal
hypersurface in $\mathbb{P}^{5}$ of degree $5$.
\end{abstract}

\maketitle

We assume that all varieties are projective, normal and defined
over $\mathbb{C}$.

\section{Introduction.}
\label{section:introduction}

The following result is proved in \cite{IsMa71}.

\begin{theorem}
\label{theorem:quartic} Every smooth hypersurface in
$\mathbb{P}^{4}$ of degree $4$ is birationally
superrigid\footnote{Let $V$ be a Fano variety such that the
singularities of the variety $V$ are at most terminal and
$\mathbb{Q}$-fac\-to\-rial singularities, and the equality
$\mathrm{rk}\,\mathrm{Pic}(V)=1$ holds. Then $V$ is called
bi\-ra\-ti\-onally rigid if it cannot be fibred into uniruled
varieties by a non-trivial rational map, and $V$ is not birational
to a Fano variety with terminal $\mathbb{Q}$-factorial
singularities of Picard rank $1$ not biregular to $V$. The variety
$V$ is called bi\-ra\-ti\-onally superrigid if it is birationally
rigid and every birational automorphism of $V$ is biregular.}.
\end{theorem}

The following generalization of Theorem~\ref{theorem:quartic} is
proved in \cite{Pu87}.

\begin{theorem}
\label{theorem:Pukhlikov} Every smooth hypersurface in
$\mathbb{P}^{5}$ of degree $5$ is birationally superrigid.
\end{theorem}

The following generalization of Theorem~\ref{theorem:quartic} is
proved in \cite{Pu88b} and \cite{Me03}.

\begin{theorem}
\label{theorem:Mella} Every $\mathbb{Q}$-factorial
nodal\,\footnote{A variety is called nodal if it has at most
isolated ordinary double points.} hypersurface in $\mathbb{P}^{4}$
of degree $4$ is birationally rigid.
\end{theorem}

In this paper we prove the following result.

\begin{theorem}
\label{theorem:quintic} Every nodal hypersurface in
$\mathbb{P}^{5}$ of degree $5$ is birationally superrigid.
\end{theorem}

It must be pointed out, that the proof of
Theorem~\ref{theorem:quintic} is based on the Shokurov
con\-nec\-ted\-ness principle (see \cite{Sho93} and Theorems~7.4
and 7.5 in \cite{Ko97}) and the Corti
inequality~(see~The\-orem~3.1 in \cite{Co00}). The proof of
Theorem~\ref{theorem:quintic} and technique used in \cite{Ch00b}
imply that every nodal quintic fourfold can not be birationally
transformed into an elliptic fibration.

\section{Smooth points.}
\label{section:smooth-points}

Let $X$ be a fourfold, $\mathcal{M}$ be a linear system on $X$
that does not have fixed components, and $O$ be a smooth point of
$X$ such that $O$ is a center of canonical singularities of the
movable log pair $(X, \lambda\mathcal{M})$, but the singularities
of the log pair $(X, \lambda\mathcal{M})$ are log terminal in a
punctured neighborhood of the point $O$, where $\lambda$ is a
positive rational number.

Let $\pi:V\to X$ be a blow up of the point $O$, and $E$ be the
$\pi$-exceptional divisor. Then
$$
K_{V}+\lambda\mathcal{B}\sim_{\mathbb{Q}}\pi^{*}\Big(K_{X}+\lambda\mathcal{M}\Big)+\big(3-m\big)E,
$$
where $\mathcal{B}$ is a proper transform of $\mathcal{M}$ on the
fourfold $V$, and $m$ is a positive rational number such that
$m/\lambda$ is the multiplicity of a general divisor of
$\mathcal{M}$ in the point $O$. Then $m>1$.

\begin{remark}
\label{remark:stupid-inequality} In fact, the inequality
$\mathrm{mult}_{O}(M_{1}\cdot M_{2})\geqslant 4/\lambda^{2}$ holds
(see Corollary~3.4 in \cite{Co00}), where $M_{1}$ and $M_{2}$ are
sufficiently general divisors of the linear system $\mathcal{M}$.
\end{remark}

In this section we prove the following result.

\begin{theorem}
\label{theorem:smooth-points} There is a line $L\subset
E\cong\mathbb{P}^{3}$  such that
$$
\mathrm{mult}_{O}\Big(M_{1}\cdot M_{2}\cdot Y\Big)\geqslant \frac{8}{\lambda^{2}},%
$$
where $M_{1}$ and $M_{2}$ are sufficiently general divisors in the
linear system $\mathcal{M}$, and $Y$ is an effective divisor on
$X$ such that
$\mathrm{dim}(\mathrm{Supp}(Y)\cap\mathrm{Supp}(M_{1}\cdot
M_{2}))=1$ and $L\subset\mathrm{Supp}(\breve{Y})$, where
$\breve{Y}$ is the proper transform of the threefold $Y$ on the
fourfold $V$.
\end{theorem}

The claim of Theorem~\ref{theorem:smooth-points} is obvious when
$m\geqslant 3$. Moreover, it follows from The\-o\-rem~7.4 and the
proof of Corollary~3.5 in \cite{Co00} that one of the following
possibilities holds:
\begin{itemize}
\item the inequality $m\geqslant 3$ holds;%
\item there is a surface $S\subset E$ such that $S$ is a center of
log canonical singularities of the log pair $(V,
\lambda\mathcal{B}+(m-2)E)$;%
\item there is a line $L\subset E\cong\mathbb{P}^{3}$ such that
$L$ is a center of log canonical singularities of the log pair
$(V, \lambda\mathcal{B}+(m-2)E)$.
\end{itemize}

\begin{lemma}
\label{lemma:smooth-points-I} Suppose that there is a surface
$S\subset E$ such that $S$ is a center of log canonical
singularities of the log pair $(V, \lambda\mathcal{B}+(m-2)E)$.
Then
$$
\mathrm{mult}_{O}\Big(M_{1}\cdot M_{2}\Big)\geqslant \frac{8}{\lambda^{2}},%
$$
where $M_{1}$ and $M_{2}$ are general divisors in the linear
system $\mathcal{M}$.
\end{lemma}

\begin{proof}
Let $B_{i}$ be a proper transform of the divisor $M_{i}$ on the
fourfold $V$. Then
$$
\mathrm{mult}_{S}\Big(B_{1}\cdot B_{2}\Big)\geqslant\frac{4\big(3-m\big)}{\lambda^{2}}%
$$
by Theorem~3.1 in \cite{Co00}. Therefore, we have
$$
\mathrm{mult}_{O}\Big(M_{1}\cdot
M_{2}\Big)\geqslant\mathrm{mult}^{2}_{O}\big(M_{i}\big)+\mathrm{mult}_{S}\Big(B_{1}\cdot
B_{2}\Big)\geqslant \frac{m^{2}+4\big(3-m\big)}{\lambda^{2}}\geqslant \frac{8}{\lambda^{2}},%
$$
which concludes the proof.
\end{proof}

Now we suppose that $m<3$ and there are no two-dimensional centers
of log canonical singularities of the log pair $(V,
\lambda\mathcal{B}+(m-2)E)$ that are contained in the
$\pi$-exceptional divisor $E$. Therefore, there is a line
$L\subset E\cong\mathbb{P}^{3}$ such that $L$ is a center of log
canonical singularities of the log pair $(V,
\lambda\mathcal{B}+(m-2)E)$.

Let $\eta:W\to V$ be a blow up of the curve $L$, and $F$ be the
$\eta$-exceptional divisor. Then
$$
K_{W}+\lambda\mathcal{D}+\big(m-3\big)\bar{E}+\big(m+n-5\big)F\sim_{\mathbb{Q}}\big(\pi\circ\eta\big)^{*}\Big(K_{X}+\lambda\mathcal{M}\Big)
$$
where $\mathcal{D}$ and $\bar{E}$ are proper transforms of the
linear system $\mathcal{M}$ and the $\pi$-exceptional divisor $E$
on the fourfold $W$ respectively, and $n$ is a positive rational
number such that the number $n/\lambda$ is the multiplicity of a
general divisor of the linear system $\mathcal{B}$ in a general
point of the curve $L$. Therefore, we have
$$
K_{W}+\lambda\mathcal{D}+\bar{H}+\big(m-2\big)\bar{E}+\big(m+n-4\big)F\sim_{\mathbb{Q}}\big(\pi\circ\eta\big)^{*}\Big(K_{X}+\lambda\mathcal{M}+H\Big),
$$
where $H$ is a sufficiently general hyperplane section of the
fourfold $X$ passing through the point $O$, and $\bar{H}$ is a
proper transform of $H$ on the fourfold $W$. Moreover, we have
$$
K_{W}+\lambda\mathcal{D}+\bar{Y}+\big(m-2\big)\bar{E}+\big(m+n-3\big)F\sim_{\mathbb{Q}}\big(\pi\circ\eta\big)^{*}\Big(K_{X}+\lambda\mathcal{M}+Y\Big),
$$
where $Y$ is a general hyperplane section of $X$ such that $L$ is
contained in the proper transform of $Y$ on the fourfold $V$, and
$\bar{Y}$ is a proper transform of $Y$ on the fourfold $W$.

\begin{lemma}
\label{lemma:smooth-points-II} Suppose that either $m+n\geqslant
4$ or there is a surface $S\subset F$ that is a center of log
canonical singularities of $(W,
\lambda\mathcal{D}+(m-2)\bar{E}+(m+n-3)F)$ and $\eta(S)=L$. Then
$$
\mathrm{mult}_{O}\Big(M_{1}\cdot M_{2}\cdot Y\Big)\geqslant \frac{8}{\lambda^{2}},%
$$
where $M_{1}$ and $M_{2}$ are general divisors in the linear
system $\mathcal{M}$.
\end{lemma}

\begin{proof}
Let $\breve{Y}$ be a proper transform of the divisor $Y$ on the
fourfold $V$. Then
$$
K_{V}+\lambda\mathcal{B}+\breve{Y}+\big(m-2\big)E\sim_{\mathbb{Q}}\pi^{*}\Big(K_{X}+\lambda\mathcal{M}+Y\Big)
$$
and $L$ is a center of log canonical singularities of the log pair
$(V, \lambda\mathcal{B}+\breve{Y}+(m-2)E)$.

The morphism $\pi\vert_{\breve{Y}}:\breve{Y}\to Y$ is a blow up of
the point $O$, the linear system $\mathcal{B}\vert_{\breve{Y}}$
does not have fixed components due to the generality in the choice
of the divisor $Y$, and we can identify the divisor $G$ with the
$\pi\vert_{\breve{Y}}$-exceptional divisor. Let us show that $L$
is a center of log canonical singularities of the log pair
$(\breve{Y},
\lambda\mathcal{B}\vert_{\breve{Y}}+(m-2)E\vert_{\breve{Y}})$.

In the case when $m+n\geqslant 4$, the equivalence
$$
K_{\bar{Y}}+\lambda\mathcal{D}\vert_{\bar{Y}}+\big(m-2\big)\bar{E}\vert_{\bar{Y}}\sim_{\mathbb{Q}}\big(\eta\vert_{\bar{Y}}\big)^{*}\Big(K_{\breve{Y}}+\lambda\mathcal{B}\vert_{\breve{Y}}+(m-2)E\vert_{\breve{Y}}\Big)+\big(3-m-n\big)F\vert_{\bar{Y}}
$$
shows that $L$ is a center of log canonical singularities of
$(\breve{Y},
\lambda\mathcal{B}\vert_{\breve{Y}}+(m-2)E\vert_{\breve{Y}})$.

Suppose that there is a surface $S\subset F$ that is a center of
log canonical singularities of the log pair $(W,
(m-2)\bar{E}+(m+n-3)F)$ and $\eta(S)=L$. Then every irreducible
component of $S\cap\bar{Y}$ is a center of log canonical
singularities of the log pair
$$
\Big(\bar{Y},
\lambda\mathcal{D}\vert_{\bar{Y}}+\big(m-2\big)\bar{E}\vert_{\bar{Y}}+\big(3-m-n\big)F\vert_{\bar{Y}}\Big),
$$
which implies that $L$ is a center of log canonical singularities
of $(\breve{Y},
\lambda\mathcal{B}\vert_{\breve{Y}}+(m-2)E\vert_{\breve{Y}})$,
because every irreducible component of the intersection
$S\cap\bar{Y}$ dominates the curve $L$.

Let $B_{i}$ be a proper transform of the divisor $M_{i}$ on the
fourfold $V$. Then
$$
\mathrm{mult}_{L}\Big(B_{1}\vert_{\breve{Y}}\cdot B_{2}\vert_{\breve{Y}}\Big)\geqslant\frac{4\big(3-m\big)}{\lambda^{2}}%
$$
by Theorem~3.1 in \cite{Co00}. Therefore, we have
$$
\mathrm{mult}_{O}\Big(M_{1}\cdot M_{2}\cdot
Y\Big)=\mathrm{mult}_{O}\Big(M_{1}\vert_{Y}\cdot
M_{2}\vert_{Y}\Big)\geqslant\frac{m^{2}}{\lambda^{2}}+\mathrm{mult}_{L}\Big(B_{1}\vert_{\breve{Y}}\cdot
B_{2}\vert_{\breve{Y}}\Big)\geqslant\frac{8}{\lambda^{2}},%
$$
which concludes the proof.
\end{proof}

Suppose that $m+n<4$.  Then in order to prove
Theorem~\ref{theorem:smooth-points} we must show that there is a
surface $S\subset F$ such that $Z$ is a center of log canonical
singularities of the log pair
$$
\Big(W,
\lambda\mathcal{D}+\big(m-2\big)\bar{E}+\big(m+n-3\big)F\Big)
$$
and $\eta(S)=L$. However, the last assertion is local and we may
assume that $X\cong\mathbb{C}^{4}$.

The singularities of the log pair $(H,
\lambda\mathcal{M}\vert_{H})$ are log terminal in a punctured
neighborhood of the point $O$. Moreover, the point $O$ is a center
of log canonical singularities of the log pair $(H,
\lambda\mathcal{M}\vert_{H})$ by Theorem~7.5 in \cite{Ko97}.
Therefore, the equivalence
$$
K_{\bar{H}}+\lambda\mathcal{D}\vert_{H}+\big(m-2\big)\bar{E}\vert_{\bar{H}}+\big(m+n-4\big)F\vert_{\bar{H}}\sim_{\mathbb{Q}}\big(\pi\circ\eta\vert_{\bar{H}}\big)^{*}\Big(K_{H}+\lambda\mathcal{M}\vert_{H}\Big),
$$
and Theorem~7.4 in \cite{Ko97} imply that we have the following
possibilities:
\begin{itemize}
\item there is a curve $C\subset F\vert_{\bar{H}}$ such that $C$
is a center of log canonical singularities of the log pair
$(\bar{H},
\lambda\mathcal{D}\vert_{H}+(m-2)\bar{E}\vert_{\bar{H}}+(m+n-4)F\vert_{\bar{H}})$;%
\item there is a point $P\in F\cap \bar{H}$ such that $P$ is a
center of log canonical singularities of the log pair $(\bar{H},
\lambda\mathcal{D}\vert_{H}+(m-2)\bar{E}\vert_{\bar{H}}+(m+n-4)F\vert_{\bar{H}})$,
and there are no other centers of log canonical singularities of
the log pair $(\bar{H},
\lambda\mathcal{D}\vert_{H}+(m-2)\bar{E}\vert_{\bar{H}}+(m+n-4)F\vert_{\bar{H}})$
except the point $P$ that are contained in the intersection $F\cap
\bar{H}$.
\end{itemize}

\begin{remark}
\label{remark:smooth-points-simple-remark} In the case when there
is a curve $C\subset F\cap \bar{H}$ such that $C$ is a center of
log canonical singularities of $(\bar{H},
\lambda\mathcal{D}\vert_{H}+(m-2)\bar{E}\vert_{\bar{H}}+(m+n-4)F\vert_{\bar{H}})$,
the curve $C$ is an intersection of the divisor $\bar{H}$ with a
surface $S\subset F$ such that $S$ is a center of log canonical
singularities of $(W, \lambda\mathcal{D}+(m-2)\bar{E}+(m+n-4)F)$
and $\eta(S)=L$.
\end{remark}

To prove Theorem~\ref{theorem:smooth-points} we may assume that
there is a point $P\in F\cap\bar{H}$ that is a center of log
canonical singularities of the log pair
$$
\Big(\bar{H},
\lambda\mathcal{D}\vert_{H}+\big(m-2\big)\bar{E}\vert_{\bar{H}}+\big(m+n-4\big)F\vert_{\bar{H}}\Big),
$$
but singularities of the log pair $(\bar{H},
\lambda\mathcal{D}\vert_{H}+(m-2)\bar{E}\vert_{\bar{H}}+(m+n-4)F\vert_{\bar{H}})$
are log terminal in a punctured neighborhood of the point
$P\in\bar{H}$.

\begin{remark}
\label{remark:smooth-points-existence-of-a-section} The morphism
$\eta\vert_{F}:F\to L$ is a $\mathbb{P}^{2}$-bundle, and the
intersection  $\bar{H}\cap F$ is just a fiber of $\eta\vert_{F}$.
Hence, the generality of in the choice of the divisor $H$ implies
the existence of a curve $Z\subset F$ such that $P=Z\cap\bar{H}$
and $Z$ is a center of log canonical singularities of the log pair
$(W, \lambda\mathcal{D}+(m-2)\bar{E}+(m+n-4)F)$.
\end{remark}

The curve $Z$ is a section of the $\mathbb{P}^{2}$-bundle
$\eta\vert_{F}$ and a center of log canonical singularities of the
log pair $(W, \lambda\mathcal{D}+(m-2)\bar{E}+(m+n-3)F)$. It
follows from Lemma~\ref{lemma:smooth-points-II} that we may assume
that $F$ does not contains surfaces dominating $L$ that are
centers of log canonical singularities of $(W,
\lambda\mathcal{D}+(m-2)\bar{E}+(m+n-3)F)$.

It follows from Theorem~7.5 in \cite{Ko97} that the point $O$ is
an isolated center of log canonical singularities of the log pair
$(Y, \lambda\mathcal{M}\vert_{Y})$. Therefore, the equivalence
$$
K_{\bar{Y}}+\lambda\mathcal{D}\vert_{Y}+\big(m-2\big)\bar{E}\vert_{\bar{Y}}+\big(m+n-3\big)F\vert_{\bar{Y}}\sim_{\mathbb{Q}}\big(\pi\circ\eta\vert_{\bar{Y}}\big)^{*}\Big(K_{Y}+\lambda\mathcal{M}\vert_{Y}\Big),
$$
and Theorem~7.4 in \cite{Ko97} imply that we have the following
possibilities:
\begin{itemize}
\item the curve $Z$ is contained in the threefold $\bar{Y}$ and
$Z$ is a center of log canonical singularities of the log pair
$(\bar{Y},
\lambda\mathcal{D}\vert_{Y}+(m-2)\bar{E}\vert_{\bar{Y}}+(m+n-3)F\vert_{\bar{Y}})$;%
\item the intersection $Z\cap \bar{Y}$ consists of a single point
that is a center of log canonical singularities of the log pair
$(\bar{Y},
\lambda\mathcal{D}\vert_{Y}+(m-2)\bar{E}\vert_{\bar{Y}}+(m+n-3)F\vert_{\bar{Y}})$.
\end{itemize}

\begin{corollary}
\label{corollary:smooth-points-main-corollary} Either $Z\subset
\bar{Y}$, or the intersection $Z\cap\bar{Y}$ consists of a single
point.
\end{corollary}

By construction we have $L\cong Z\cong\mathbb{P}^{1}$ and
$$
F\cong\mathrm{Proj}\Big(\mathcal{O}_{\mathbb{P}^{1}}(-1)\oplus\mathcal{O}_{\mathbb{P}^{1}}(1)\oplus\mathcal{O}_{\mathbb{P}^{1}}(1)\Big),%
$$
but the equivalence $\bar{Y}\vert_{F}\sim B+D$ holds, where $B$ is
the tautological line bundle on the threefold $F$, and $D$ is a
fiber of the natural projection $\eta\vert_{F}:F\to
L\cong\mathbb{P}^{1}$.

\begin{lemma}
\label{lemma:smooth-points-III} The equality
$h^{1}(\mathcal{O}_{W}(\bar{Y}-F))$ holds.
\end{lemma}

\begin{proof}
The divisor $-\eta^{*}(E)-F$ intersects every curve contained in
$\bar{E}$ non-negatively and
$$
\Big(-\eta^{*}\big(E\big)-F\Big)\vert_{F}\sim B+D,
$$
which implies that $-4\eta^{*}(E)-4F$ is $(\pi\circ\eta)$-big and
$(\pi\circ\eta)$-nef. However, we have
$$
K_{W}-4\big(\pi\circ\eta\big)^{*}\big(E\big)-4F\sim \bar{Y}-F
$$
and $X\cong\mathbb{C}^{4}$, which implies that
$h^{1}(\mathcal{O}_{W}(\bar{Y}-F))=0$ by the Kawamata--Vieh\-weg
vanishing theorem (see Theorem~2.3 in \cite{Ko97}).
\end{proof}

Thus, the restriction map $H^{0}(\mathcal{O}_{W}(\bar{Y}))\to
H^{0}(\mathcal{O}_{F}(\bar{Y}\vert_{F}))$ is surjective, but the
complete linear system $\vert \bar{Y}\vert_{F}\vert$ does not have
base points.

\begin{corollary}
\label{corollary:smooth-points-freeness} The intersection
$\bar{Y}\cap Z$ consists of a single point.
\end{corollary}

Let $\mathcal{I}_{Z}$ be an ideal sheaf of $Z$ on $F$. Then
$R^{1}\,(\eta\vert_{F})_{*}(B\otimes\mathcal{I}_{Z})=0$ and there
is a surjective map
$\psi:\mathcal{O}_{\mathbb{P}^{1}}(-1)\oplus\mathcal{O}_{\mathbb{P}^{1}}(1)\oplus\mathcal{O}_{\mathbb{P}^{1}}(1)\to\mathcal{O}_{\mathbb{P}^{1}}(k)$,
where $k=B\cdot Z$. The map $\psi$ is given by a an element of the
group
$$
H^{0}\Big(\mathcal{O}_{\mathbb{P}^{1}}\big(k+1\big)\Big)\oplus H^{0}\Big(\mathcal{O}_{\mathbb{P}^{1}}\big(k-1\big)\Big)\oplus H^{0}\Big(\mathcal{O}_{\mathbb{P}^{1}}\big(k-1\big)\Big),%
$$
which implies that $k\geqslant -1$.

\begin{lemma}
\label{lemma:smooth-points-IV} The equality $k=0$ is impossible.
\end{lemma}

\begin{proof}
Suppose $k=0$. Then $\psi$ is given by matrix $(ax+by, 0, 0)$,
where $a$ and $b$ are constants and $(x:y)$ are homogeneous
coordinates on $L\cong\mathbb{P}^{1}$. Thus $\psi$ is not
surjective over the point of $L$ at which $ax+by$ vanishes.
\end{proof}

Therefore, the divisor $B$ can not have trivial intersection with
$Z$. Hence the intersection of the divisor $\bar{Y}$ with the
curve $Z$ is either trivial or consists of more than one point,
but the intersection $\bar{Y}\cap Z$ consists of one point. The
obtained contradiction proves Theorem~\ref{theorem:smooth-points}.

\section{Singular points.}
\label{section:singular-points}

Let $X$ be a fourfold, $\mathcal{M}$ be a linear system on the
fourfold $X$ that does not have fixed components, and $O$ be an
isolated ordinary double point of $X$ such that $O$ is a center of
canonical singularities of the log pair $(X, \lambda\mathcal{M})$,
but $(X, \lambda\mathcal{M})$ has log terminal singularities in a
punctured neighborhood of the point $O$, where $\lambda$ is a
positive rational number.

Let $\pi:V\to X$ be a blow up of the point $O$, and $E$ be the
exceptional divisor of the birational morphism $\pi$. Then $E$ is
a smooth quadric hypersurface in $\mathbb{P}^{4}$ and
$$
K_{V}+\lambda\mathcal{B}\sim_{\mathbb{Q}}\pi^{*}\Big(K_{X}+\lambda\mathcal{M}\Big)+\big(2-m\big)E,
$$
where $\mathcal{B}$ is a proper transform of $\mathcal{M}$ on the
fourfold $V$, and $m$ is a po\-si\-tive rational number. It
follows from Theorem~3.10 in \cite{Co00} and Theorem~7.5 in
\cite{Ko97} that $m>1$.

Let $M_{1}$ and $M_{2}$  be general divisors in $\mathcal{M}$, and
$Y$ be a hyperplane section of $X$. Then
$$
B_{i}\sim\pi^{*}\big(M_{i}\big)-\frac{m}{\lambda}E
$$
and $\breve{Y}\sim\pi^{*}(Y)-E$, where $B_{i}$ and $\breve{Y}$ are
proper transforms of the divisors $M_{i}$ and $Y$ on the fourfold
$V$ respectively. Suppose in addition that the following
conditions hold:
\begin{itemize}
\item the point $O$ is an ordinary double point of the threefold $Y$;%
\item the threefold $Y$ does not contain surfaces contained in the base locus of $\mathcal{M}$;%
\item the threefold $\breve{Y}$ does not contain surfaces contained in the base locus of $\mathcal{B}$.%
\end{itemize}

Let $H$ be a sufficiently general hyperplane section of the
fourfold $X$ that passes through the point $O$, and
$\Sigma=\mathrm{Supp}(B_{1}\cdot B_{2}\cdot \breve{Y}\cdot
\breve{H})\cap\mathrm{Supp}(E)$, where $\breve{H}$ is a proper
transform of the divisor $H$ on the fourfold $V$. Then
$|\Sigma|<+\infty$. Put
$$
\nabla=2\frac{m^{2}}{\lambda^{2}}+\sum_{P\in \Sigma}\mathrm{mult}_{P}\Big(B_{1}\cdot B_{2}\cdot\breve{Y}\cdot\breve{H}\Big).%
$$

In this section we prove the following result.

\begin{theorem}
\label{theorem:singular-points} There is a line $L\subset
E\subset\mathbb{P}^{4}$ such that $\nabla>6n^{2}$ whenever
$L\subset\mathrm{Supp}(\breve{Y})$.%
\end{theorem}

Let $H^{\prime}$ be a sufficiently general hyperplane section of
the fourfold $X$ that passes through the point $O$, and
$\Sigma^{\prime}=\mathrm{Supp}(B_{1}\cdot B_{2}\cdot
\breve{H}^{\prime}\cdot \breve{H})\cap\mathrm{Supp}(E)$, where
$\breve{H}^{\prime}$ is a proper transform of the divisor
$H^{\prime}$ on the fourfold $V$. Then
$|\Sigma^{\prime}|<+\infty$. Put
$$
\nabla^{\prime}=2\frac{m^{2}}{\lambda^{2}}+\sum_{P\in \Sigma^{\prime}}\mathrm{mult}_{P}\Big(B_{1}\cdot B_{2}\cdot\breve{H}^{\prime}\cdot\breve{H}\Big),%
$$
which implies that the inequality $\nabla\geqslant\nabla^{\prime}$
holds.

In order to prove Theorem~\ref{theorem:singular-points} we may
assume that $m<2$. Then the singularities of the log pair $(V,
\lambda\mathcal{B}+(m-1)E)$ are not log terminal in the
neighborhood of $E$.

\begin{lemma}
\label{lemma:singular-points-I} Suppose that there is a surface
$S\subset E$ such that $S$ is a center of log canonical
singularities of the log pair $(V, \lambda\mathcal{B}+(m-1)E)$.
Then $\nabla\geqslant 6/\lambda^{2}$.
\end{lemma}

\begin{proof}
It follows from Theorem~3.1 in \cite{Co00} that the inequality
$$
\mathrm{mult}_{S}\Big(B_{1}\cdot B_{2}\Big)\geqslant\frac{4\big(2-m\big)}{\lambda^{2}}%
$$
holds. Therefore, we have
$$
\nabla^{\prime}=2\frac{m^{2}}{\lambda^{2}}+\sum_{P\in \Sigma^{\prime}}\mathrm{mult}_{P}\Big(B_{1}\cdot B_{2}\cdot\breve{H}^{\prime}\cdot\breve{H}\Big)\geqslant 2\frac{m^{2}}{\lambda^{2}}+\mathrm{mult}_{S}\Big(B_{1}\cdot B_{2}\Big)\geqslant\frac{2m^{2}+4\big(2-m\big)}{\lambda^{2}}\geqslant 6/\lambda^{2},%
$$
which concludes the proof.
\end{proof}

Therefore, in order to prove Theorem~\ref{theorem:singular-points}
we may assume that the set of centers of log canonical
singularities of the log pair $(V, \lambda\mathcal{B}+(m-1)E)$
does not contains surfaces that are contained in $E$. Then the
claim of Theorem~7.4 in \cite{Ko97} together with the equivalences
$$
K_{V}+\breve{H}^{\prime}+\lambda\mathcal{B}+\big(m-1\big)E\sim_{\mathbb{Q}}\pi^{*}\Big(K_{X}+H^{\prime}+\lambda\mathcal{M}\Big)
$$
and
$K_{\breve{H}^{\prime}}+\lambda\mathcal{B}\vert_{\breve{H}^{\prime}}+(m-1)E\vert_{\breve{H}^{\prime}}\sim_{\mathbb{Q}}\pi^{*}(K_{H^{\prime}}+\lambda\mathcal{M}\vert_{\breve{H}^{\prime}})$
imply that there is a line $L\subset E\subset\mathbb{P}^{4}$ such
that $L$ is the unique center of log canonical singularities of
$(V, \lambda\mathcal{B}+(m-1)E)$ that is contained in $E$, because
$H^{\prime}$ is sufficiently general, but the point $O$ is a
center of log canonical singularities of the log pair
$(H^{\prime}, \lambda\mathcal{M}\vert_{\breve{H}^{\prime}})$ by
Theorem~7.5 in \cite{Ko97}.

Now we suppose that $L\subset\mathrm{Supp}(\breve{Y})$. Then $L$
is a center of log canonical singularities of the log pair
$(\breve{Y},
\lambda\mathcal{B}\vert_{\breve{Y}}+(m-1)E\vert_{\breve{Y}})$ by
Theorem~7.5 in \cite{Ko97}. Hence, we have
$$
\mathrm{mult}_{L}\Big(B_{1}\cdot B_{2}\cdot\breve{Y}\Big)=\mathrm{mult}_{L}\Big(B_{1}\vert_{\breve{Y}}\cdot B_{2}\vert_{\breve{Y}}\Big)\geqslant\frac{4\big(2-m\big)}{\lambda^{2}}%
$$
by Theorem~3.1 in \cite{Co00}. Therefore, we have
$$
\nabla=2\frac{m^{2}}{\lambda^{2}}+\sum_{P\in\Sigma}\mathrm{mult}_{P}\Big(B_{1}\cdot B_{2}\cdot\breve{Y}\cdot\breve{H}\Big)\geqslant 2\frac{m^{2}}{\lambda^{2}}+\mathrm{mult}_{L}\Big(B_{1}\cdot B_{2}\cdot\breve{Y}\Big)\geqslant\frac{2m^{2}+4\big(2-m\big)}{\lambda^{2}}\geqslant 6/\lambda^{2},%
$$
which concludes the proof of
Theorem~\ref{theorem:singular-points}.

\section{Birational rigidity.}
\label{section:birational-rigidity}

In this section we prove Theorem~\ref{theorem:quintic}. Let $X$ be
a hypersurface in $\mathbb{P}^{5}$ of degree $5$ with at most
isolated ordinary double points. Then the group $\mathrm{Cl}(X)$
is generated by the class of a hyperplane section (see
\cite{CalLy94}). Suppose that the quintic fourfold $X$ is not
birationally super\-rigid. Then there is a linear system
$\mathcal{M}$ on the fourfold $X$ that does not have fixed
components, but the singularities of the log pair $(X,
\frac{1}{n}\mathcal{M})$ are not canonical (see \cite{Co95}),
where $n$ is a natural number such that the equivalence
$\mathcal{M}\sim -nK_{X}$ holds.

Let $Z$ be an irreducible subvariety of the fourfold $X$ having
maximal dimension such that the singularities of the log pair $(X,
\frac{1}{n}\mathcal{M})$ are not canonical in a general point of
the subvariety $Z$. Then $\mathrm{mult}_{Z}(\mathcal{M})>n$, which
implies that $\mathrm{dim}(Z)\leqslant 1$ due to \cite{Pu98a}.

\begin{lemma}
\label{lemma:birational-rigidity-smooth-points} The subvariety $Z$
is not a smooth point of the hypersurface $X$.
\end{lemma}

\begin{proof}
Suppose that $Z$ is a smooth point of the hypersurface $X$. Let
$\pi:V\to X$ be a blow up of the point $Z$, and $E$ be the
exceptional divisor of the morphism $\pi$. Then
$$
K_{V}+{\frac{1}{n}}\mathcal{B}\sim_{\mathbb{Q}} \pi^{*}\Big(K_{X}+{\frac{1}{n}}\mathcal{M}\Big)+\Big(3-\mathrm{mult}_{Z}\big(\mathcal{M}\big)/n\Big)E\sim_{\mathbb{Q}}\Big(3-\mathrm{mult}_{Z}\big(\mathcal{M}\big)/n\Big)E,%
$$
where $\mathcal{B}$ is a proper transform of the linear system
$\mathcal{M}$ on the variety $V$.

Let $M_{1}$ and $M_{2}$ be general divisors in $\mathcal{M}$, and
$H_{1}$ and $H_{2}$ be general hyperplane sections of the
hypersurface $X$ passing through the point $Z$. Then
$$
5n^{2}=M_{1}\cdot M_{2}\cdot H_{1}\cdot H_{2}\geqslant \mathrm{mult}_{Z}\Big(M_{1}\cdot M_{2}\Big)\geqslant\mathrm{mult}^{2}_{Z}\big(\mathcal{M}\big)>n^2,%
$$
which implies that $\mathrm{mult}_{Z}(\mathcal{M})\leqslant
\sqrt{5}n<3n$.

Now it follows from Theorem~\ref{theorem:smooth-points} that there
is a line $L\subset E\cong\mathbb{P}^{3}$ such that
$$
\mathrm{mult}_{Z}\Big(M_{1}\cdot M_{2}\cdot Y\Big)>8n^{2},
$$
where $Y$ is a hyperplane section of the hypersurface $X$ such
that
$$
\mathrm{dim}\Big(\mathrm{Supp}\big(Y\big)\cap\mathrm{Supp}\Big(M_{1}\cdot M_{2}\Big)\Big)=1%
$$
and $L\subset\mathrm{Supp}(\breve{Y})$, where $\breve{Y}$ is the
proper transform of $Y$ on the fourfold $V$.

Let $\mathcal{D}$ be a linear subsystem in
$|\mathcal{O}_{\mathbb{P}^{5}}(1)\vert_{X}|$ such that
$$
D\in \mathcal{D}\iff L\subset\mathrm{Supp}\big(\breve{D}\big)\
\mathrm{or}\ \mathrm{mult}_{Z}\big(D\big)\geqslant 2,
$$
where $\breve{D}$ is a proper transform of $D$ on the fourfold
$V$. Then there is a two-dimensional linear subspace
$\Pi\subset\mathbb{P}^{5}$ such that the base locus of
$\mathcal{D}$ consists of $X\cap\Pi$.

Suppose that $\Pi\not\subset\mathrm{Supp}(M_{1}\cdot M_{2})$. Let
$D$ be a general divisor in $\mathcal{D}$. Then
$$
5n^{2}=M_{1}\cdot M_{2}\cdot D\cdot H_{1}\geqslant \mathrm{mult}_{Z}\Big(M_{1}\cdot M_{2}\cdot D\Big)>8n^{2},%
$$
which is a contradiction. In particular, the quintic $X$ contains
the plane $\Pi$.

Let $\bar{X}$ be a general hyperplane section of $X$ containing
$\Pi$. Then $\bar{X}$ is a quintic hyper\-sur\-face with isolated
singularities in $\mathbb{P}^{4}$ that is smooth at $Z$. Let
$\bar{\pi}:\bar{V}\to\bar{X}$ be a blow up of the point $Z$, and
$\bar{E}$ be the $\bar{\pi}$-exceptional divisor. There is a
commutative diagram
$$
\xymatrix{
&\bar{V}\ar@{->}[d]_{\bar{\pi}}\ar@{^{(}->}[rr]&& V\ar@{->}[d]^{\pi}&\\%
&\bar{X}\ar@{^{(}->}[rr]&&X,&}
$$
where we identify $\bar{V}$ with the proper transform of $\bar{X}$
on the fourfold $V$, and $\bar{E}=E\cap\bar{V}$.

The plane $\Pi$ is a fixed component of the linear system
$\mathcal{M}\vert_{\bar{X}}$. Moreover, we have
$$
\mathcal{M}\vert_{\bar{X}}=\mathcal{R}+\alpha\Pi,
$$
where $\mathcal{R}$ is a linear system on $\bar{X}$ that does not
have fixed components, and $\alpha$ is a multiplicity of a general
divisor of the linear system $\mathcal{M}$ in a general point of
the plane $\Pi$.

Let $\mathcal{L}$ and $\breve{\Pi}$ be the proper transforms of
the linear system $\mathcal{R}$ and the plane $\Pi$ on the
threefold $\bar{V}$ respectively. Then $L\subset\breve{\Pi}$ and
$$
K_{\bar{V}}+\frac{1}{n}\Big(\mathcal{L}+\alpha\breve{\Pi}\Big)\sim_{\mathbb{Q}}
\bar{\pi}^{*}\Big(K_{\bar{X}}+\frac{1}{n}\mathcal{R}+\frac{\alpha}{n}\Pi\Big)+
\Big(2-\mathrm{mult}_{Z}\big(\mathcal{R}\big)/n-\alpha/n\Big)\bar{E},
$$
but the proof of Theorem~\ref{theorem:smooth-points} implies the
singularities of the log pair
$$
\Big(\bar{V},
\frac{1}{n}\mathcal{L}+\frac{\alpha}{n}\bar{\Pi}+\frac{\mathrm{mult}_{Z}\big(\mathcal{R}\big)-\alpha-2n}{n}\bar{E}\Big)
$$
in a general point of the curve $L$. Therefore, the inequality
$$
\mathrm{mult}_{L}\Big(L_{1}\cdot L_{2}\Big)>4\Big(3n-\mathrm{mult}_{Z}\big(\mathcal{R}\big)-\alpha\Big)\big(n-\alpha\big)%
$$
holds by Theorem~3.1 in \cite{Co00}, where $L_{1}$ and $L_{2}$ are
general surfaces in $\mathcal{L}$. We have
$$
\mathrm{mult}_{Z}\Big(R_{1}\cdot R_{2}\Big)\geqslant\mathrm{mult}^{2}_{Z}\big(\mathcal{R}\big)+\mathrm{mult}_{L}\Big(L_{1}\cdot L_{2}\Big)>\mathrm{mult}^{2}_{Z}\big(\mathcal{R}\big)+4\Big(3n-\mathrm{mult}_{Z}\big(\mathcal{R}\big)-\alpha\Big)\big(n-\alpha\big),%
$$
where $R_{1}=\bar{\pi}(L_{1})$ and $R_{2}=\bar{\pi}(L_{2})$.

Let $S$ be a general hyperplane section of $\bar{Y}$ that contains
$Z$. Then
$$
5n^{2}-2n\alpha-3\alpha^{2}=R_{1}\vert_{S}\cdot R_{2}\vert_{S}\geqslant \mathrm{mult}_{Z}\Big(R_{1}\cdot R_{2}\Big)>\Big(\mathrm{mult}^{2}_{Z}\big(\mathcal{R}\big)+\alpha-n\Big)^{2}+8\big(n-\alpha\big),%
$$
which is a contradiction.
\end{proof}

\begin{lemma}
\label{lemma:birational-rigidity-singular-points} The subvariety
$Z$ is not a singular point of the hypersurface $X$.
\end{lemma}

\begin{proof}
Suppose that $Z$ is a singular point of the quintic $X$. Let
$\pi:V\to X$ be a blow up of the point $Z$, and $E$ be the
exceptional divisor of $\pi$. Then $E$ is a quadric in
$\mathbb{P}^{4}$ and
$$
K_{V}+\frac{1}{n}\mathcal{B}\sim_{\mathbb{Q}}\pi^{*}\Big(K_{X}+\frac{1}{n}\mathcal{M}\Big)+\big(2-m\big)E,
$$
where $\mathcal{B}$ is a proper transform of $\mathcal{M}$ on the
fourfold $V$, and $m$ is a positive ra\-ti\-onal number. The
inequality $m>1$ holds by Theorem~3.10 in \cite{Co00} and
Theorem~7.5 in \cite{Ko97}.

Let $M_{1}$ and $M_{2}$ be general divisors of the linear system
$\mathcal{M}$, and $B_{i}$  be a proper transform of the divisor
$M_{i}$ on the fourfold $V$. Then
$$
B_{1}\cdot B_{2}\cdot\breve{H}_{1}\cdot\breve{H}_{2}\leqslant 5n^{2}-2m^{2}n^{2},%
$$
where $\breve{H}_{1}$ and $\breve{H}_{2}$ are proper transforms on
the fourfold $V$ of hyperplane sections of the hypersurface $X$
that pass through the point $Z$ such that the equality
$$
\mathrm{dim}\Big(\mathrm{Supp}\Big(B_{1}\cdot
B_{2}\Big)\cap\mathrm{Supp}\big(\breve{H}_{1}\big)\cap\mathrm{Supp}\big(\breve{H}_{2}\big)\Big)=0
$$
holds. In particular, the inequality $m\leqslant\sqrt{5/2}$ holds.

It follows from Theorem~\ref{theorem:singular-points} that there
is a line $L\subset E\subset\mathbb{P}^{4}$ such that
$$
B_{1}\cdot B_{2}\cdot\breve{Y}\cdot\breve{H}>6n^{2}-2m^{2}n^{2},%
$$
where $\breve{H}$ is a proper transform on $V$  of a general
hyperplane sections of $X$ that passes through the point $Z$, and
$\breve{Y}$ is a proper transform on $V$  of a hyperplane sections
of the fourfold $X$ passing through $Z$ such that $Z$ is an
ordinary double point of $\pi(\breve{Y})$, the equality
$$
\mathrm{dim}\Big(\mathrm{Supp}\Big(B_{1}\cdot B_{2}\Big)\cap\mathrm{Supp}\big(\breve{Y}\big)\Big)=1%
$$
holds and $L\subset\mathrm{Supp}(\breve{Y})$.

Let $\mathcal{D}$ be a linear subsystem in
$|\mathcal{O}_{\mathbb{P}^{5}}(1)\vert_{X}|$ spanned by the
divisors whose proper transforms on $V$ contain $L$. Then there is
a two-dimensional linear subspace $\Pi\subset\mathbb{P}^{5}$ such
that the base locus of $\mathcal{D}$ consists of $X\cap\Pi$.
Therefore, we have $\Pi\subset\mathrm{Supp}(M_{1}\cdot M_{2})$.

Let $\bar{X}$ be a general hyperplane section of $X$ that contains
$\Pi$. Then $\bar{X}$ is a quintic hyper\-sur\-face in
$\mathbb{P}^{4}$ having isolated singularities, and $Z$ is an
isolated ordinary double point of the quintic $\bar{X}$, which has
at most canonical singularities (see Corollary~4.9 in
\cite{Ko97}).

The plane $\Pi$ is a fixed component of the linear system
$\mathcal{M}\vert_{\bar{X}}$. Moreover, we have
$$
\mathcal{M}\vert_{\bar{X}}=\mathcal{R}+\alpha\Pi,
$$
where $\mathcal{R}$ is a linear system on $\bar{X}$ that does not
have fixed components, and $\alpha$ is a multiplicity of a
sufficiently general divisor of the linear system $\mathcal{M}$ in
a general point of the plane $\Pi$. The plane $\Pi$ and a general
surface of $\mathcal{R}$ are not $\mathbb{Q}$-Cartier divisors on
$\bar{X}$.

Let $\bar{\pi}:\bar{V}\to\bar{X}$ be a composition of the blow up
of $Z$ with a sub\-sequent
$\mathbb{Q}$-fac\-to\-ri\-a\-li\-zation, and $\bar{E}$ be the
$\bar{\pi}$-exceptional divisor. Then $\bar{E}$ is a smooth
quadric surface and
$$
K_{\bar{V}}+\frac{1}{n}\mathcal{L}+\frac{\alpha}{n}\breve{\Pi}\sim_{\mathbb{Q}}
\bar{\pi}^{*}\Big(K_{\bar{X}}+\frac{1}{n}\mathcal{R}+\frac{\alpha}{n}\Pi\Big)+
\big(1-m\big)\bar{E},
$$
where $\mathcal{L}$ and $\breve{\Pi}$ are proper transforms of
$\mathcal{R}$ and $\Pi$ on the threefold $\bar{V}$ respectively.

Let $\bar{L}=\bar{E}\cap\bar{\Pi}$. Then the curve $\bar{L}$ is a
line on the quadric $\bar{E}$. Moreover, it follows from the proof
of Theorem~\ref{theorem:singular-points}  that the singularities
of the log pair
$$
\Big(\bar{V},
\frac{1}{n}\mathcal{L}+\frac{\alpha}{n}\bar{\Pi}+\big(m-1\big)\bar{E}\Big)
$$
are not log canonical in a general point of $\bar{L}$. We have
$\mathrm{mult}_{\bar{L}}(\mathcal{L})>2n-\alpha-mn$, but
$$
\mathrm{mult}_{\bar{L}}\Big(L_{1}\cdot
L_{2}\Big)>4\big(2n-mn\big)\big(n-\alpha\big)%
$$
by Theorem~3.1 in \cite{Co00}, where $L_{1}$ and $L_{2}$ are
general surfaces in $\mathcal{L}$. The inequality
$$
3\alpha-mn+n>\mathrm{mult}_{\bar{L}}\big(\mathcal{L}\big)\geqslant 0%
$$
holds, because $\bar{C}\cdot L_{i}\geqslant
\mathrm{mult}_{\bar{L}}(\mathcal{L})$ and
$\bar{C}\cdot\bar{\Pi}=-3$, where $\bar{C}$ is a proper transform
on the threefold $\bar{V}$ of a sufficiently general line
contained in $\Pi$. Thus, we have $\alpha>n/4$.

Let $\psi:\bar{X}\dasharrow\mathbb{P}^{1}$ be a projection from
the plane $\Pi$. Then $\psi$ is not defined in the points where
the plane $\Pi$ is not a Cartier divisor on $\bar{X}$. In
particular, the rational map $\psi$ is not defined in the point
$Z$. However, we may assume that the birational morphism
$\bar{\pi}$ resolves the indeterminacy of the rational map $\psi$.
Therefore, we have a commutative diagram
$$
\xymatrix{
&&\bar{V}\ar@{->}[dl]_{\bar{\pi}}\ar@{->}[dr]^{\eta}&&\\%
&\bar{X}\ar@{-->}[rr]_{\psi}&&\mathbb{P}^{2},&}
$$
where $\eta$ is a morphism. Let $S$ be a sufficiently general
fiber of $\eta$ and $H$ be a hyperplane section of $\bar{X}$. Then
$\bar{\pi}(S)$ is a quartic surface in $\mathbb{P}^{4}$ and $S\sim
\bar{\pi}^{*}(H)-\bar{E}-\bar{\Pi}$, which implies
\begin{equation}
\label{equation:main-inequality}
\Big(\bar{\pi}^{*}\big(H\big)-\bar{E}-\bar{\Pi}\Big)\cdot\Big(\bar{\pi}^{*}\big(nH\big)-nm\bar{E}-\alpha\bar{\Pi}\Big)^{2}>4\big(2n-mn\big)\big(n-\alpha\big),%
\end{equation}
because $S\cdot\bar{L}=1$ and
$L_{i}\sim_{\mathbb{Q}}\bar{\pi}^{*}(nH)-nm\bar{E}-\alpha\bar{\Pi}$.

We have $H^{3}=5$, $\bar{E}^{3}=2$,
$\bar{\pi}^{*}(H)\cdot\bar{\Pi}^{2}=-3$,
$\bar{E}\cdot\bar{\Pi}^{2}=0$, $\bar{\Pi}\cdot\bar{E}^{2}=-1$, and
$$
\bar{\Pi}^{3}=\Big(\bar{\pi}^{*}\big(H\big)-\bar{E}-S\Big)\cdot\Pi^{2}=-3-S\cdot\bar{\Pi}^{2},
$$
but $\bar{\pi}(S)\cap\Pi$ is a hyperplane section of the surface
$\bar{\pi}(S)\subset\mathbb{P}^{4}$. Moreover, the generality in
the choice of the threefold $\bar{X}$ implies that the threefold
$\bar{X}$ has isolated ordinary double points, the quartic surface
$\bar{\pi}(S)$ is smooth, and the birational morphism
$\bar{\pi}\vert_{S}$ is a blow up of a point $Z$ on the surface
$\bar{\pi}(S)$. Therefore, we have $S\cdot\bar{\Pi}^{2}=3$, which
gives $\bar{\Pi}^{3}=-6$.

The inequality~\ref{equation:main-inequality} implies that the
inequality
$$
4n^{2}-9\alpha^{2}+4n\alpha-2mn\alpha-m^{2}n^{2}>
4\big(2n-mn\big)\big(n-\alpha\big)
$$
holds. Therefore, we have
$$
0>4n^{2}+9\alpha^{2}-12n\alpha+6mn\alpha-4mn^{2}+m^{2}n^{2}=\Big(3\alpha-2n+mn\Big)^{2},
$$
which is a contradiction.
\end{proof}

Therefore, the subvariety $Z$ is an irreducible curve.

\begin{lemma}
\label{lemma:birational-rigidity-plane-curve} The curve $Z$ is a
line.
\end{lemma}

\begin{proof}
Suppose that $Z$ is not a line. Let $P_{1}$ and $P_{2}$ be
sufficiently general points of the curve $Z$, and $L$ be a line in
$\mathbb{P}^{5}$ that passes through the points $P_{1}$ and
$P_{2}$. Then $L\ne Z$.

Let $H_{1}$ and $H_{2}$ be sufficiently general hyperplane section
of the fourfold $X$ that pass through $P_{1}$ and $P_{2}$. Put
$S=H_{1}\cap H_{2}$. Then the singularities of the log pair $(S,
\frac{1}{n}\mathcal{M}\vert_{S})$ are not log canonical in the
points $P_{1}$ and $P_{2}$ by Theorem~7.5 in \cite{Ko97}, but the
singularities of the log pair $(S,
\frac{1}{n}\mathcal{M}\vert_{S})$ are log canonical in punctured
neighborhoods of these points, because the secant variety of the
curve $Z$ is at least two-dimensional. Put
$$
\mathcal{M}\vert_{S}=\mathcal{B}+\gamma L,
$$
where $\mathcal{B}$ is a linear system on $S$ that does not have
fixed components, and $\gamma$ is the multiplicity of a general
divisor of $\mathcal{M}$ in a general point of the line $L$. Then
$$
\mathrm{mult}_{P_{i}}\Big(B_{1}\cdot B_{2}\Big)>4\big(n^{2}-\gamma n\big)%
$$
by Theorem~3.1 in \cite{Co00}, where $B_{1}$ and $B_{2}$ are
general divisors in $\mathcal{B}$. Thus, we have
$$
5n^{2}-2\gamma n-3\gamma^{2}=B_{1}\cdot
B_{2}\geqslant\mathrm{mult}_{P_{1}}\Big(B_{1}\cdot
B_{2}\Big)+\mathrm{mult}_{P_{2}}\Big(B_{1}\cdot
B_{2}\Big)>8\big(n^{2}-\gamma n\big),
$$
which is a contradiction.
\end{proof}

Let $Y$ be a sufficiently general hyperplane section of the
fourfold $X$ that passes through the line $Z$, and
$\mathcal{B}=\mathcal{M}\vert_{Y}$. Then $Y$ is a quintic
threefold in $\mathbb{P}^{4}$, the linear system $\mathcal{B}$
does not have fixed components, but the singularities of the log
pair $(Y, \frac{1}{n}\mathcal{B})$ are not log canonical in a
general point of the curve $Z$ by Theorem~7.5 in \cite{Ko97}.

The line $Z$ contains a singular point of the threefold $Y$ due to
\cite{Pu98a}, but $Z$ contains at most $4$ singular points of the
threefold $Y$. Put
$$
Z\cap \mathrm{Sing}\big(Y\big)=\big\{P_{1},\ldots , P_{k}\big\},
$$
where $P_{i}$ is a singular point of the threefold $Y$ and
$k\leqslant 4$. Then the point $P_{i}$ is an isolated ordinary
double point of the threefold $Y$, and the group $\mathrm{Cl}(Y)$
is generated by the class of a hyperplane section (see
\cite{Cy01}).

Let $\pi:V\to Y$ be a blow up of the points $\{P_{1},\ldots ,
P_{k}\}$, and $E_{i}$ be an ex\-cep\-ti\-onal divisor of the
morphism $\pi$ such that $\pi(E_{i})=P_{i}$. Then
$$
D\sim\pi^{*}\Big(\mathcal{O}_{\mathbb{P}^{4}}\big(n\big)\vert_{Y}\Big)-\sum_{i=1}^{4}m_{i}E_{i},
$$
where $D$ is a proper transform on $V$ of a general divisor in
$\mathcal{B}$, and $m_{i}$ is a natural number.

Let $\bar{Z}$ be a proper transform of $Z$ on the fourfold $V$.
Then $V$ is smooth and
$$
2m_{i}\geqslant \mathrm{mult}_{Q_{i}}\big(D\big)\geqslant \mathrm{mult}_{\bar{Z}}\big(D\big)=\mathrm{mult}_{Z}\big(\mathcal{B}\big)\geqslant n,%
$$
where $Q_{i}=\bar{Z}\cap E_{i}$. Hence, we have
$m_{i}\geqslant\mathrm{mult}_{Z}(\mathcal{B})/2\geqslant n/2$.

Let $\Pi$ be a general plane in $\mathbb{P}^{4}$ that contains the
line $Z$, and $C$ be a quartic curve in the plane $\Pi$ such that
$\Pi\cap Y=L\cup C$. Then $|C\cap Z|=4$ and the curve $C$ contains
all singular points of the threefold $Y$ that is contained in $Z$.
Let $\bar{C}$ be a proper transform of the curve $C$ on the
threefold $V$. Then $\bar{C}\not\subset\mathrm{Supp}(D)$ and
$$
0\leqslant D\cdot\bar{C}=4n-\sum_{i=1}^{k}m_{i}\leqslant
4n-2k\mathrm{mult}_{Z}\big(\mathcal{B}\big)+\big(4-k\big)\mathrm{mult}_{Z}\big(\mathcal{B}\big),
$$
which implies that $\mathrm{mult}_{Z}(\mathcal{B})\leqslant 2n$.

Let $\omega:W\to V$ be a blow up of $\bar{Z}$, and $G$ be the
$\omega$-exceptional divisor. Then
$$
K_{W}+\frac{1}{n}\mathcal{D}\sim
\big(\pi\circ\omega\big)^{*}\Big(K_{Y}+\frac{1}{n}\mathcal{B}\Big)+\sum_{i=1}^{k}\Big(1-\frac{m_{i}}{n}\Big)\breve{E}_{i}+\Big(1-\mathrm{mult}_{Z}\big(\mathcal{B}\big)/n\Big)G,
$$
where $\breve{E}_{i}$ and $\mathcal{D}$ are proper transforms of
the divisor $E_{i}$ and the linear system $\mathcal{B}$ on the
threefold $W$ respectively. Therefore, the inequality
$\mathrm{mult}_{Z}(\mathcal{B})<2n$ implies the existence of an
irreducible curve $L\subset G$ such that $\omega(L)=\bar{Z}$, but
the sin\-gu\-la\-rities of the log pair
$$
\Big(W,\frac{1}{n}\mathcal{D}+\frac{\mathrm{mult}_{Z}\big(\mathcal{B}\big)-n}{n}G\Big)
$$
are not log canonical in a general point of $L$. Thus, we have
$\mathrm{mult}_{L}(\mathcal{D})+\mathrm{mult}_{Z}(\mathcal{B})>2n$.

The surface $E_{i}$ is isomorphic to
$\mathbb{P}^{1}\times\mathbb{P}^{1}$. Let $A_{i}$ and $B_{i}$ be
the fibers of the projections of the surface $E_{i}$ to
$\mathbb{P}^{1}$ that pass through the point $Q_{i}$, and
$\bar{A}_{i}$ and $\bar{B}_{i}$ be proper transforms of the curves
$A_{i}$ and $B_{i}$ on the threefold $W$ respectively. Then
$$
\mathcal{N}_{W/\bar{A}_{i}}\cong\mathcal{N}_{W/\bar{B}_{i}}\cong\mathcal{O}_{\mathbb{P}^{1}}(-1)\oplus\mathcal{O}_{\mathbb{P}^{1}}(-1),
$$
which implies that we can flop the curves $\bar{A}_{i}$ and
$\bar{B}_{i}$. Namely, let $\xi:U\to W$ be a blow up of the curves
$\bar{A}_{1}, \bar{B}_{1}, \ldots, \bar{A}_{k}, \bar{B}_{k}$, and
$F_{i}$ and $H_{i}$ be the exceptional divisors of the morphism
$\xi$ such that $\xi(F_{i})=\bar{A}_{i}$ and
$\xi(H_{i})=\bar{B}_{i}$. Then $F_{i}\cong
H_{i}\cong\mathbb{P}^{1}\times\mathbb{P}^{1}$ and there is a
birational morphism $\xi^{\prime}:U\to W^{\prime}$ such that
$\xi^{\prime}(F_{i})$ and $\xi^{\prime}(H_{i})$ are rational
curves, but the map $\xi^{\prime}\circ\xi^{-1}$ is not an
isomorphism in the neighborhood of the curves $\bar{A}_{i}$ and
$\bar{B}_{i}$.

Let $E_{i}^{\prime}$ be a proper transform of $E_{i}$ on the
threefold $W^{\prime}$. Then $E_{i}^{\prime}\cong\mathbb{P}^{2}$
and we can contract the surface $E_{i}^{\prime}$ to a singular
point of type $\frac{1}{2}(1,1,1)$, because
$$
\mathcal{N}_{W^{\prime}/E_{i}^{\prime}}\cong\mathcal{O}_{E_{i}^{\prime}}\Big(E_{i}^{\prime}\vert_{E_{i}^{\prime}}\Big)\cong\mathcal{O}_{\mathbb{P}^{2}}(-2).
$$

Let $\omega^{\prime}:W^{\prime}\to V^{\prime}$ be a
con\-trac\-tion of $E_{1}^{\prime},\ldots, E_{k}^{\prime}$, and
$G^{\prime}$ be a proper trans\-form of the surface $G$ on the
threefold $V^{\prime}$. Then there is a birational morphism
$\pi^{\prime}:V^{\prime}\to Y$ that contracts the divisor
$G^{\prime}$ to the line $Z$. Hence, we constructed the
commutative diagram
$$
\xymatrix{
&&W\ar@{->}[dl]_{\omega}&&U\ar@{->}[ll]_{\xi}\ar@{->}[dr]^{\xi^{\prime}}\\%
&V\ar@{->}[dr]_{\pi}&&&&W^{\prime}\ar@{->}[dl]^{\omega^{\prime}}\\%
&&Y&&V^{\prime}\ar@{->}[ll]_{\pi^{\prime}}&}
$$
such that $V^{\prime}$ is projective and $\mathbb{Q}$-factorial,
and $\mathrm{rk}\,\mathrm{Pic}(V^{\prime})=2$. Therefore, the
birational morphism $\pi^{\prime}:V^{\prime}\to Y$ is an exremal
terminal divisorial contraction (see \cite{Tz03}). We have
$$
K_{V^{\prime}}+\frac{1}{n}\mathcal{R}\sim_{\mathbb{Q}}
{\pi^{\prime}}^{*}\Big(K_{Y}+\frac{1}{n}\mathcal{B}\Big)+\Big(1-\mathrm{mult}_{Z}\big(\mathcal{B}\big)/n\Big)G^{\prime},
$$
where $\mathcal{R}$ is a proper transform of $\mathcal{B}$ on the
threefold $V^{\prime}$. Let $L^{\prime}$ be a proper transform of
the curve $L$ on the threefold $V^{\prime}$. Then
$(V^{\prime},\frac{1}{n}\mathcal{R}+(\mathrm{mult}_{Z}(\mathcal{B})/n-1)G^{\prime})$
is not log canonical in a general point of the curve $L^{\prime}$.
Hence, the inequality
$$
\mathrm{mult}_{L^{\prime}}\Big(R_{1}\cdot R_{2}\Big)>4n\Big(2n-\mathrm{mult}_{Z}\big(\mathcal{B}\big)\Big)%
$$
holds by Theorem~3.1 in \cite{Co00}, where $R_{1}$ and $R_{2}$ are
general surfaces in $\mathcal{R}$.

Let $H$ be a hyperplane section of the threefold $Y$, and
$P^{\prime}_{i}=\pi^{\prime}(E_{i}^{\prime})$. Then the base locus
of the linear system $|{\pi^{\prime}}^{*}(H)-G^{\prime}|$ consists
of the points $P^{\prime}_{1},\ldots, P^{\prime}_{k}$. The
construction of the morphism $\pi^{\prime}$ implies that
${G^{\prime}}^{3}=2-k/2$. We have
$({\pi^{\prime}}^{*}(H)-G^{\prime})\cdot L^{\prime}=0$, because
$$
2\Big({\pi^{\prime}}^{*}(H)-G^{\prime}\Big)\cdot R_{1}\cdot
R_{2}=2\Big({\pi^{\prime}}^{*}(H)-G^{\prime}\Big)\cdot
\Big({\pi^{\prime}}^{*}(nH)-\mathrm{mult}_{Z}\big(\mathcal{B}\big)G^{\prime}\Big)^{2}<
4n\Big(2n-\mathrm{mult}_{Z}\big(\mathcal{B}\big)\Big)
$$
and $2({\pi^{\prime}}^{*}(H)-G^{\prime})$ is a Cartier divisor. In
particular, we have $L^{\prime}\cap\{P^{\prime}_{1},\ldots,
P^{\prime}_{k}\}=\emptyset$.

\begin{corollary}
\label{corollary:birational-rigidity-Mori-cone} The cone
$\overline{\mathbb{NE}}(V^{\prime})$ is generated by the curves
$\xi^{\prime}(F_{i})\equiv\xi^{\prime}(H_{i})$ and $L^{\prime}$.
\end{corollary}

Let $\bar{H}$ be a proper transform on the threefold $V$ of a
sufficiently general hyperplane section of the threefold $Y$ that
contains the line $Z$. Then $\bar{H}$ is a smooth surface such
that the equality $\bar{Z}^{2}=-3$ holds on the surface $\bar{H}$.
One the other hand, we have
$$
c_{1}\Big(\mathcal{N}_{V/\bar{Z}}\Big)=-2-K_{V}\cdot\bar{Z}=-2-k,
$$
which implies that $\mathcal{N}_{V/\bar{Z}}\cong
\mathcal{O}_{\mathbb{P}^{1}}(a)\oplus\mathcal{O}_{\mathbb{P}^{1}}(b)$,
where $a$ and $b$ are integer numbers such that the inequality
$a\geqslant b\geqslant -3$ holds and $a+b=-2-k$.

Let $\breve{H}$ be the proper transform of the divisor $\bar{H}$
on the threefold $W$. Then
$$
\breve{H}\sim\big(\pi\circ\omega\big)^{*}\Big(\mathcal{O}_{\mathbb{P}^{4}}\big(1\big)\vert_{Y}\Big)-\sum_{i=1}^{4}\bar{E}_{i}-G
$$
and $G^{3}=2+k$. Elementary calculations implies that $\breve{H}$
intersects the curve $L$ in its general point in the case when $L$
is not the exceptional section of the projection
$$
\omega\vert_{G}:G\cong\mathbb{F}_{a-b}\to\bar{Z}
$$
or $b\ne-3$. Thus, we have $b=-3$ and the curve $L$ is the
ex\-cep\-ti\-onal section of the ruled surface
$G\cong\mathbb{F}_{4-k}$. Moreover, the construction of the map
$\omega^{\prime}\circ\xi^{\prime}\circ\xi^{-1}$ implies that
$$
L\cap\Big(\bar{A}_{i}\cup\bar{B}_{i}\Big)\ne\emptyset
$$
for every $i=1,\ldots,k$, because otherwise $P^{\prime}_{i}\in
L^{\prime}$ and $({\pi^{\prime}}^{*}(H)-G^{\prime})\cdot
L^{\prime}\ne 0$.

\begin{lemma}
\label{lemma:birational-rigidity-k-is-not-4} The inequality
$k\leqslant 3$ holds.
\end{lemma}

\begin{proof}
Suppose that $k=4$. Then $|{\pi^{\prime}}^{*}(H)-2G^{\prime}|$
contains a divisor $T$ such that $\pi^{\prime}(T)$ is a hyperplane
section of $Y$ that tangents $Y$ along $Z$. The cycle $T\cdot
R_{i}$ must be effective, but
$$
T\cdot R_{i}\equiv \Big(6n-6\mathrm{mult}_{Z}\big(\mathcal{B}\big)\Big)\xi^{\prime}\big(F_{i}\big)+\Big(5n-2\mathrm{mult}_{Z}\big(\mathcal{B}\big)\Big)L^{\prime},%
$$
which implies that $\mathrm{mult}_{Z}(\mathcal{B})\leqslant n$.
Contradiction.
\end{proof}

Let $H^{\prime}$ be a proper transform on $V^{\prime}$ of a
general hyperplane section of $Y$. Then
$$
\left\{\aligned
&H^{\prime}\cdot H^{\prime}\equiv 10\xi^{\prime}\big(F_{i}\big)+5L^{\prime},\\%
&H^{\prime}\cdot G^{\prime}\equiv 2\xi^{\prime}\big(F_{i}\big)\equiv\xi^{\prime}\big(F_{i}\big)+\xi^{\prime}\big(H_{i}\big),\\%
&G^{\prime}\cdot G^{\prime}\equiv \big(k-6\big)\xi^{\prime}\big(F_{i}\big)-L^{\prime},\\%
\endaligned
\right.
$$
which implies that the equivalence
$$
R_{1}\cdot R_{2}\equiv \Big(10n^{2}-4n\mathrm{mult}_{Z}\big(\mathcal{B}\big)+\big(k-6\big)\mathrm{mult}^{2}_{Z}\big(\mathcal{B}\big)\Big)\xi^{\prime}\big(F_{i}\big)+\Big(5n^{2}-\mathrm{mult}^{2}_{Z}\big(\mathcal{B}\big)\Big)L^{\prime}%
$$
holds. Thus, we have
$10n-4\mathrm{mult}_{Z}(\mathcal{B})+(k-6)\mathrm{mult}^{2}_{Z}(\mathcal{B})\geqslant
0$, which implies that
$$
n<\mathrm{mult}_{Z}\big(\mathcal{B}\big)\leqslant\frac{\sqrt{34}-2}{3}n<\frac{32}{25}n<\frac{4}{3}n.%
$$

Let $\pi^{\prime\prime}:V^{\prime\prime}\to V^{\prime}$ be a blow
up of $L^{\prime}$, and $G^{\prime\prime}$ be the
$\pi^{\prime\prime}$-exceptional divisor. Then
$$
K_{V^{\prime\prime}}+\frac{1}{n}\mathcal{L}+\frac{\mathrm{mult}_{Z}\big(\mathcal{B}\big)-n}{n}\bar{G}^{\prime}+
\frac{\mathrm{mult}_{L^{\prime}}\big(\mathcal{R}\big)+\mathrm{mult}_{Z}\big(\mathcal{B}\big)-2n}{n}G^{\prime\prime}
\sim_{\mathbb{Q}}
\big({\pi^{\prime}\circ\pi^{\prime\prime}}\big)^{*}\Big(K_{Y}+\frac{1}{n}\mathcal{B}\Big),
$$
where $\mathcal{L}$ and $\bar{G}^{\prime}$ are proper transforms
of $\mathcal{R}$ and $G^{\prime}$ on the threefold
$V^{\prime\prime}$ respectively, which implies that either the
inequality
$\mathrm{mult}_{L^{\prime}}(\mathcal{R})+\mathrm{mult}_{Z}(\mathcal{B})>3n$
holds, or the log pair
\begin{equation}
\label{equation:the-log-pair} \Big(V^{\prime\prime},
\frac{1}{n}\mathcal{L}+\Big(\mathrm{mult}_{Z}\big(\mathcal{B})/n-1\big)\bar{G}^{\prime}+
\Big(\mathrm{mult}_{L^{\prime}}\big(\mathcal{R}\big)/n+\mathrm{mult}_{Z}\big(\mathcal{B}\big)/n-2\Big)G^{\prime\prime}\Big)
\end{equation}
are not log canonical in a general point of a curve dominating the
curve $L^{\prime}$. We have
$$
\mathrm{mult}_{L^{\prime}}(\mathcal{R})\leqslant\mathrm{mult}_{Z}(\mathcal{B})\leqslant\frac{\sqrt{34}-2}{3}n<\frac{4}{3}n,
$$
which implies the inequality
$\mathrm{mult}_{L^{\prime}}(\mathcal{R})+\mathrm{mult}_{Z}(\mathcal{B})<3n$.
Therefore, there is an irreducible curve $L^{\prime\prime}\subset
G^{\prime\prime}$ such that
$\pi^{\prime\prime}(L^{\prime\prime})=L^{\prime}$, but the
singularities of the log pair~\ref{equation:the-log-pair} are not
log canonical in a general point of the curve $L^{\prime\prime}$.
In particular, the inequality
$$
\mathrm{mult}_{L^{\prime\prime}}\big(\mathcal{L}\big)>4n-\mathrm{mult}_{L^{\prime}}\big(\mathcal{R}\big)-2\mathrm{mult}_{Z}\big(\mathcal{B}\big)
$$
holds, because the inequality
$\mathrm{mult}_{L^{\prime}}(\mathcal{R})>2n-\mathrm{mult}_{Z}(\mathcal{B})$
holds, but the singularities of the log
pair~\ref{equation:the-log-pair} are not canonical in a general
point of the curve $L^{\prime\prime}$.

\begin{lemma}
\label{lemma:birational-rigidity-L-prime-prime} The curve
$L^{\prime\prime}$ is not contained in the divisor
$\bar{G}^{\prime}$.
\end{lemma}

\begin{proof}
Suppose that $L^{\prime\prime}=G^{\prime\prime}\cap
\bar{G}^{\prime}$. Then taking the intersection of a general
surface in the linear system $\mathcal{L}$ and a general fiber of
the morphism
$(\pi^{\prime}\circ\pi^{\prime\prime})\vert_{\bar{G}^{\prime}}$ we
see that
$$
\mathrm{mult}_{Z}\big(\mathcal{B}\big)-\mathrm{mult}_{L^{\prime}}\big(\mathcal{R}\big)>4n-\mathrm{mult}_{L^{\prime}}\big(\mathcal{R}\big)-2\mathrm{mult}_{Z}\big(\mathcal{B}\big),
$$
which implies that
$\mathrm{mult}_{Z}(\mathcal{B})>4n/3>n(\sqrt{34}-2)/3$.
Contradiction.
\end{proof}

Let $L_{1}$ and $L_{2}$ be general surfaces in the linear system
$\mathcal{L}$, and $H^{\prime\prime}$ be a proper transform of a
general hyperplane section of the threefold $Y$ on the threefold
$V^{\prime\prime}$. Then
$$
5n^{2}-\mathrm{mult}^{2}_{L^{\prime}}\big(\mathcal{R}\big)-\mathrm{mult}^{2}_{Z}\big(\mathcal{B}\big)=
H^{\prime\prime}\cdot L_{1}\cdot L_{2}>4n\Big(3n-\mathrm{mult}_{L^{\prime}}\big(\mathcal{R}\big)-\mathrm{mult}_{Z}\big(\mathcal{B}\big)\Big)%
$$
by Theorem~3.1 in \cite{Co00}. Therefore, we have
$$
n^{2}>\Big(\mathrm{mult}_{L^{\prime}}\big(\mathcal{R}\big)-2n\Big)^{2}+\Big(\mathrm{mult}_{Z}\big(\mathcal{B}\big)-2n\Big)^{2}\geqslant
\frac{2n^{2}\Big(\sqrt{34}-8\Big)^{2}}{3}>n^{2},%
$$
which is a contradiction. The obtained contradiction concludes the
proof of Theorem~\ref{theorem:quintic}.

\end{document}